\documentclass[12pt]{article}
\usepackage{amsmath,amssymb}
\newcommand{\tur}{turbulen}
\newcommand{\rey}{Reynolds }

\renewcommand{\a}{\alpha}
\renewcommand{\b}{\beta}
\newcommand{\bs}{\bigskip}

\renewcommand{\l}{\lambda}
\renewcommand{\L}{\Lambda}

\renewcommand{\th}{\theta}

\newcommand{{\z}}{\Bbb Z}

\begin{document}
\setlength{\baselineskip}{13pt}

\begin{center}
{\large\bf The Characteristic Length Scale of the}\smallskip\\
{\large\bf Intermediate Structure in Zero-Pressure-Gradient}\smallskip\\
{\large\bf Boundary Layer Flow}

\vspace{1 truein}

{\sc G. I. Barenblatt,${}^1$ \ A. J. Chorin${}^1$ and V. M.
Prostokishin${}^2$}

\vspace{.33 truein}
${}^1$Department of Mathematics and\\
Lawrence Berkeley National Laboratory \\
University of California\\
Berkeley, California 94720

\bs\bs

${}^2$P. P. Shirshov Institute of Oceanology\\Russian Academy of
Sciences\\
  36, Nakhimov Prospect\\
  Moscow 117218 Russia
\end{center}

\bs\bs\bs\begin{quote}
\begin{center}{\bf Abstract.}\end{center}
In a turbulent boundary layer over a smooth flat plate with zero pressure
gradient, the  intermediate structure between the viscous sublayer and the
free stream consists of two layers: one adjacent to the viscous sublayer and
one adjacent to the free stream. When the level of turbulence in the free
stream is low, the boundary between the two layers is sharp and both have a
self-similar structure described by Reynolds-number-dependent scaling
(power) laws. This structure introduces two length scales: one --- the wall
region thickness --- determined by the sharp boundary between the two
intermediate layers, the second determined by the condition that the
velocity distribution in the first intermediate layer be the one common to
all wall-bounded  flows, and in particular coincide with the scaling law
previously determined for pipe flows. Using recent experimental data we
determine both these length scales and show that they are close.
Our results disagree with the classical model of the ``wake region".
\end{quote}

\newpage
\section{Introduction}

Turbulent boundary layer flow over a smooth flat plate outside a close
vicinity of the plate tip contains two unambiguous elements: The viscous
sublayer  adjacent to the plate, where the velocity gradient is large and
the viscous stress is comparable with the \rey stress, and the
statistically uniform free stream.

According to classical theory [${}^1$], the region intermediate between
these two consists of two layers with different properties. The first,
adjacent to the viscous sublayer is a universal,
Reynolds-number-independent logarithmic layer. In the second, the ``wake
region", there is a smooth transition from the universal logarithmic layer
to the free stream.

Our analysis of {\it all} available experiments [${}^{2-4}$] contradicts
this classical theory. Indeed, in the clear-cut case of a smooth plate and
low free stream turbulence, the intermediate structure does consist of two
layers. However, the boundary between them is sharp. Most important, both
layers are self-similar, substantially  Reynolds-number-dependent, and
described by different scaling laws.
It is interesting to note (see the details below) that the same
configuration of two self-similar layers with a sharp interface between them
can be seen  in all runs used in [${}^1$] for the illustration of
the wake region model.

We found it possible [${}^{2-4}$] to introduce a characteristic length
scale $\L$ so that the average velocity distribution in the first
intermediate layer coincides with  Reynolds-number-dependent scaling law
obtained previously for pipe flows, when the \rey number is chosen as
$U\L/\nu$, with
$U$ the free stream velocity and $\nu$ the fluid's kinematic viscosity. The
sharp boundary between the self-similar intermediate layers also defines a
length scale $\l$. We show, by analysis of experimental data, that these two
length scales $\l$ and $\L$ are close.

\section{Background}
In a previous paper [${}^2$] we noted that when the \tur ce level in the
free stream is small, the intermediate structure between the viscous
sublayer and the free stream consists of two self-similar layers: one
adjacent to the viscous sublayer where the average velocity profile  is
described by the scaling law
\begin{equation}\label{21}
\phi=A\eta^\a \ ,
\end{equation}
and one adjacent to the free stream where
\begin{equation}\label{22}
\phi=B\eta^\b \ .
\end{equation}
Here
\[
\phi=\frac{u}{u_*} \ , \qquad
u_*=\sqrt{\frac{\tau}{\rho}} \ ,\qquad
\eta=\frac{u_*y}{\nu} \ ,
\]
$u$ is the average velocity, $\tau$ is the shear stress at the wall, $\rho$
and $\nu$ are the fluid density and kinematic viscosity; $A,B$, $\a$ and
$\b$ are \rey number dependent constants.

Our processing of all experimental data available in the literature
[${}^{3,4}$] confirmed these observations and showed that it is always
possible to find a length scale $\L$ so that, setting $Re=U\L/\nu$, we can
represent the scaling law (\ref{21}) in the form
\begin{equation}\label{23}
\phi=\left(\frac{1}{\sqrt{3}}\ln Re +\frac 52\right)\eta^{\frac{3}{2\ln Re}}
\end{equation}
obtained by us earlier for pipe flows (see, e.g. [${}^5$]). This suggests
that {\it the structure of wall regions in all wall-bounded shear flows
at large \rey numbers is identical, if the length scale and velocity scale
are properly selected.} The natural question is, however, what is the
physical meaning of this length scale $\L$ in  boundary layer flow? This
question is of substantial importance and should be clarified for proper
understanding of  the identity of scaling laws for different
wall-bounded shear flows.

We note that the intermediate structure has another characteristic length
 scale $\l$ --- the wall-region thickness determined by
the sharp intersection $\eta=\eta_*$ of the two velocity distribution laws
$\phi=A\eta^\a$ and
$\phi=B\eta^\b$ valid in the different layers. We have
\begin{equation}\label{24}
A\eta_*^\a =B\eta_*^\b \ ,
\end{equation}
so that
\begin{equation}\label{25}
\eta_*=\big(\frac AB\big)^{\frac{1}{\b-\a}} \ ,
\end{equation}
and the wall-region thickness $\l$ is determined by the relation
\begin{equation}\label{26}
\l=\big(\frac AB\big)^{\frac{1}{\b-\a}} \ \frac{\nu}{u_*} \ .
\end{equation}
On the other hand, the characteristic length scale $\L$ is determined by the
relation
\begin{equation}\label{27}
\L \ = \ Re\frac \nu U \ =
\ \big(\frac{u_*}{U}\big)\big(\frac{\nu}{u_*}\big)Re
\end{equation}
so that the ratio of these two scales is
\begin{equation}\label{28}
\frac{\L}{\l} \ = \ \left(\frac{u_*}{U} Re\right) \frac{1}{\eta_*} \
= \ \left(\frac{u_*}{U}Re\right) \big(\frac BA\big)^{\frac{1}{\b-\a}} \ .
\end{equation}

\section{Analysis of experimental data}

We analyzed the recent data of J.M.~\"Osterlund presented on the Internet
({\tt www.kth.se/$\sim$jens/zpg/}). The data  seem to us to be
reliable, however much we disagree with their processing and interpretation
in the paper by \"Osterlund  et al [${}^6$] (see [${}^4$]). All 70 runs
presented on the Internet give the characteristic broken-line  average
velocity distribution in $\lg\eta,\lg\phi$ coordinates (see the examples in
Figure 1; all the other cases are similar), so that the possibility of
determining $A,\a,B,\b$ and $\eta_*$ accurately from these experimental
data  is unquestionable. These results are presented in Table 1 for all
\"Osterlund's experiments where $Re_{\th}=U\th/\nu > 10,000$. Here $\th$ is
the momentum thickness; the runs in the \"Osterlund's experimental data are
labelled by $Re_{\th}$.

The effective \rey number $Re$ was obtained [${}^4$] by the formula
\begin{equation}\label{31}
\ln Re=\frac 12(\ln Re_1+\ln Re_2)
\end{equation}
where $\ln Re_1$ and $\ln Re_2$ are the solutions of the equations
\begin{equation}\label{32}
\frac{1}{\sqrt{3}} \ln Re_1+\frac 52=A \ , \qquad
\frac{3}{2\ln Re_2}=\a
\end{equation}
and the values of $A$ and $\a$ were obtained by standard statistical
processing of
\"Osterlund's data. For $Re_{\th}>10,000$ the difference $\delta$ between
$\ln Re_1$ and $\ln Re_2$  does not exceed 3\%, so that they coincide within
experimental accuracy.

According to (\ref{27}) and (\ref{28})
\begin{equation}\label{33}
\lg \frac{\L}{\l} \ =  \
(\lg Re-\lg\eta_*)+\lg\frac{u_*}{U} \ .
\end{equation}
The data for $u_*$ and $U$ are presented by \"Osterlund on the Internet for
each run. In Figure 2 we present the values of $\lg(\L/\l)$ for all runs.
The mean value of $\lg(\L/\l)$ is approximately 0.2, so that the
characteristic length scale $\L$ is about 1.6 times the thickness of the
wall region.

If we take into account that $\L$ is calculated from the value of $Re$, and
that $\ln Re$, not $Re$ itself, has been determined from experiment, the
ratios
$\L/\l$ as shown in Figure 2 are close to unity.

We processed in [${}^{3,4}$] the data of 90 {\it zero-pressure-gradient}
boundary layer experiments at low free stream turbulence performed by
different authors during the last 25 years --- all the experiments of this
kind available to us. Without exception,  all runs revealed identical
configurations of the intermediate structure in the  boundary layer: two
adjacent self-similar layers separated by a sharp interface.

According to the classical model [${}^1$], the
intermediate structure consists of the (universal) logarithmic layer and a
non-self-similar ``wake region" smoothly matching the logarithmic layer.
It was natural to also process the very data used in
[${}^1$] for the justification of the wake region model with the general
procedure we used on the other data. The data  presented in Figure 21 of
[${}^1$] were scanned and replotted in
$\lg\eta,\lg\phi$ coordinates as was done for all experimental data
processed in [${}^{3,4}$]. Processing revealed the same broken-line
structure, i.e.~two adjacent self-similar layers (see Table 2 and Figure 3,
where a typical example is presented). The difference between $\ln Re_1$ and
$\ln Re_2$ determined from the wall layer data is small: this shows that the
procedure is adequate.
We conjecture that the values of $\b$ are larger than in newer experiments
because of a non-zero pressure gradient in all these runs (see [${}^1$]).
The results of our processing fail to confirm the wake region model
proposed in [${}^1$].

\section{Conclusion}

We have shown that one can find a length scale $\L$ so that, if
the \rey number $Re$ in a zero-pressure-gradient boundary layer flow is
defined by
$Re=U\L/\nu$, with $U$ the free stream velocity and $\nu$ the kinematic
viscosity, then the scaling law for the self-similar region adjacent to the
viscous sublayer coincides with the scaling law for \tur t pipe flow. Using
the recent experimental data of \"Osterlund
({\tt www.kth.se/$\sim$jens/zpg/}) we confirmed this fact and reached the
important conclusion that $\L$ is roughly equal to the wall-region thickness.
Our results are in disagreement with the classical model of the wake region
in the boundary layer [${}^1$].

\bs{\bf Acknowledgments.} This work was supported in part by the
Applied Mathematics subprogram of the U.S.~Department of Energy under
contract DE--AC03--76--SF00098, and in part by the National Science
Foudnation under grants DMS\,94--16431 and DMS\,97--32710.

\bs\bs\begin{center}{\Large\bf References}\end{center}

\begin{enumerate}

\item Coles, D. (1956). {\it J. Fluid Mech.} {\bf 1}, pp.191--226.

\item Barenblatt, G.I., Chorin, A.J., Hald, O. and Prostokishin, V.M. (1997).
 {\it Proceedings National Academy of Sciences USA}
{\bf 94}, pp.7817--7819.

\item Barenblatt, G.I., Chorin, A.J., and Prostokishin, V.M. (2000).
{\it J. Fluid Mech.}, in press.

\item Barenblatt, G.I., Chorin, A.J., and Prostokishin, V.M.
(2000). Analysis of experimental investigations of self-similar structures in
zero-pressure-gradient boundary layers at large \rey numbers.
Center for Pure \& Applied Mathematics, UC Berkeley, CPAM-777.

\item Barenblatt, G.I., Chorin, A.J., and Prostokishin, V.M. (1997).
 {\it Applied
Mechanics Reviews}, vol.~{\bf 50}, no.~7, pp.413--429.

\item \"Osterlund, J.M., Johansson, A.V.,  Nagib, H.M. and Hites, M.H.
(2000).  {\it Physics of Fluids}, vol.~{\bf 12}, no.~1, pp.1--4.
\end{enumerate}

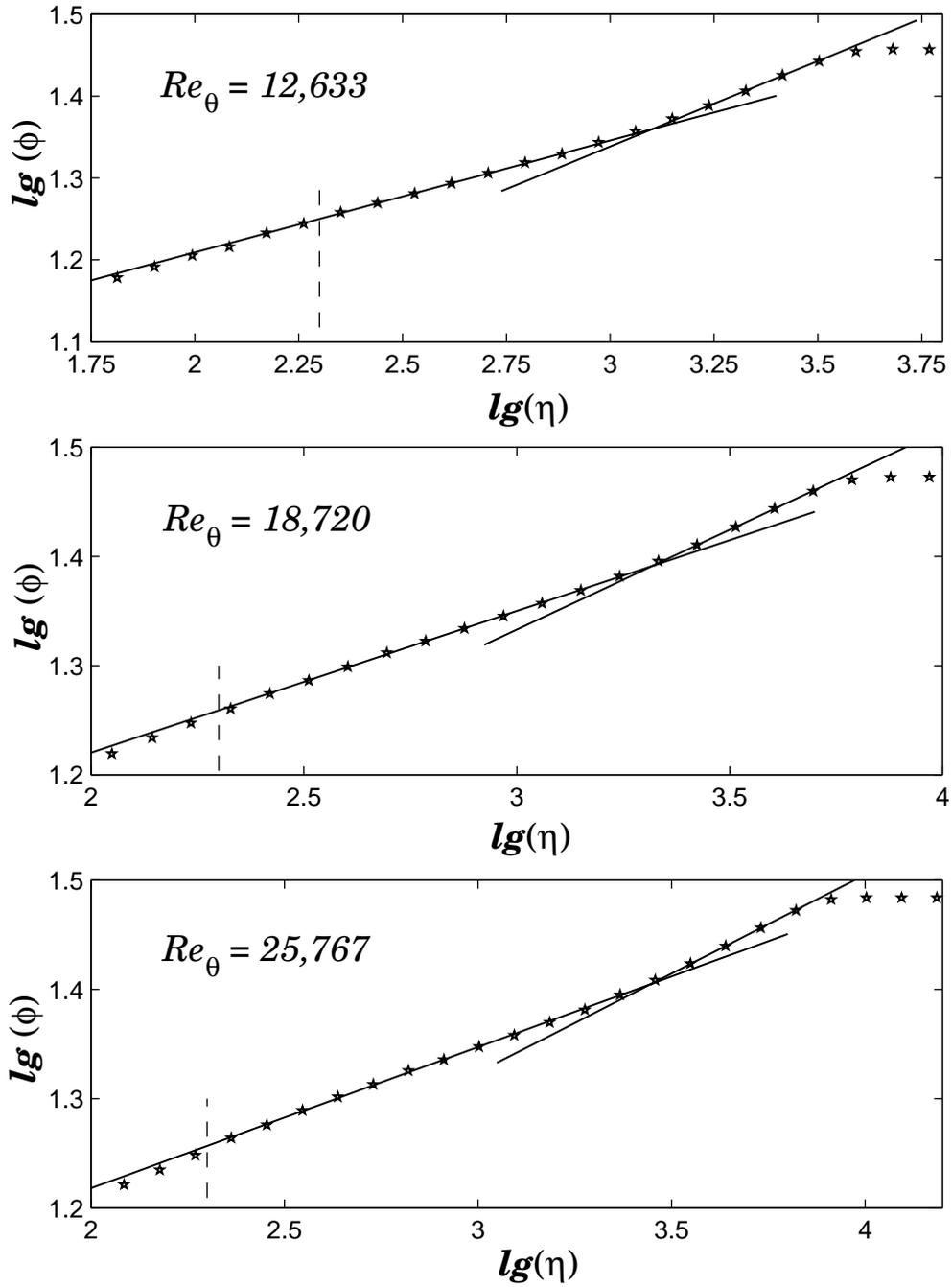
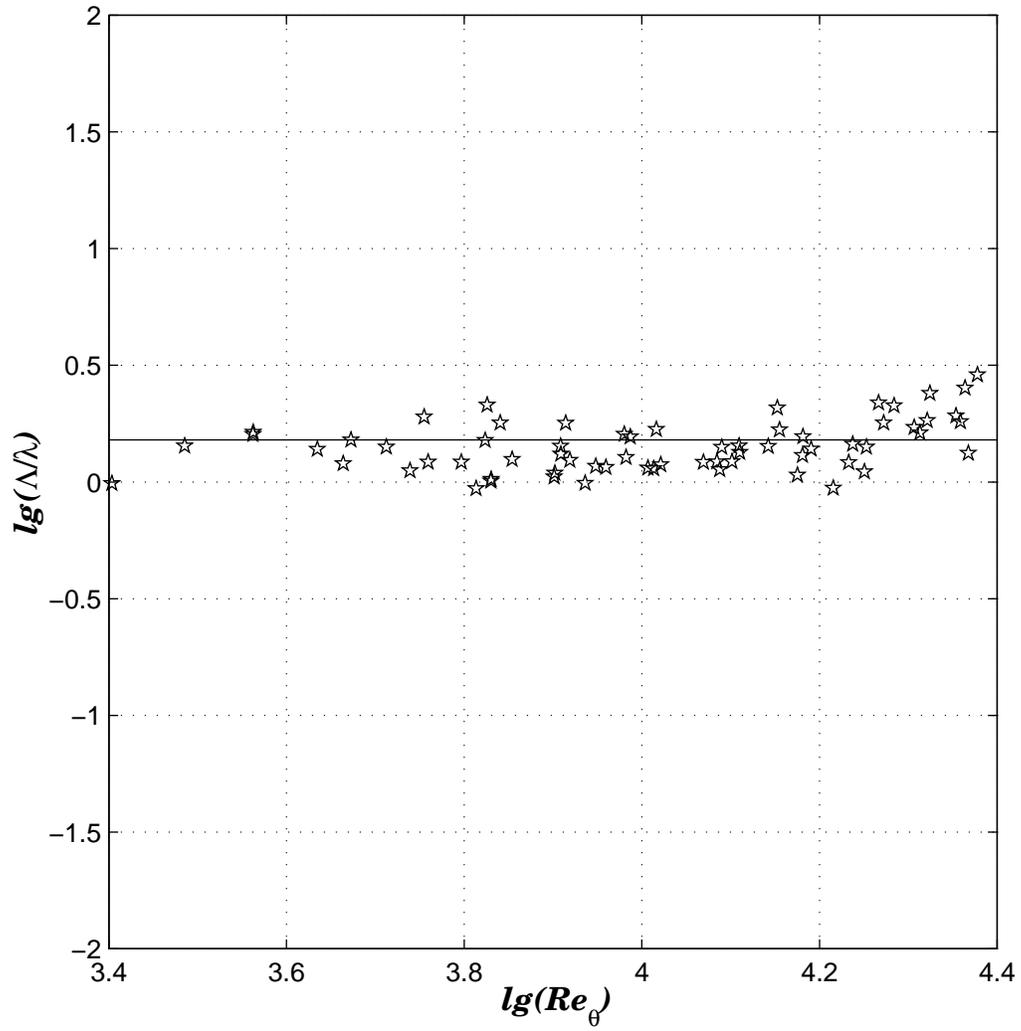
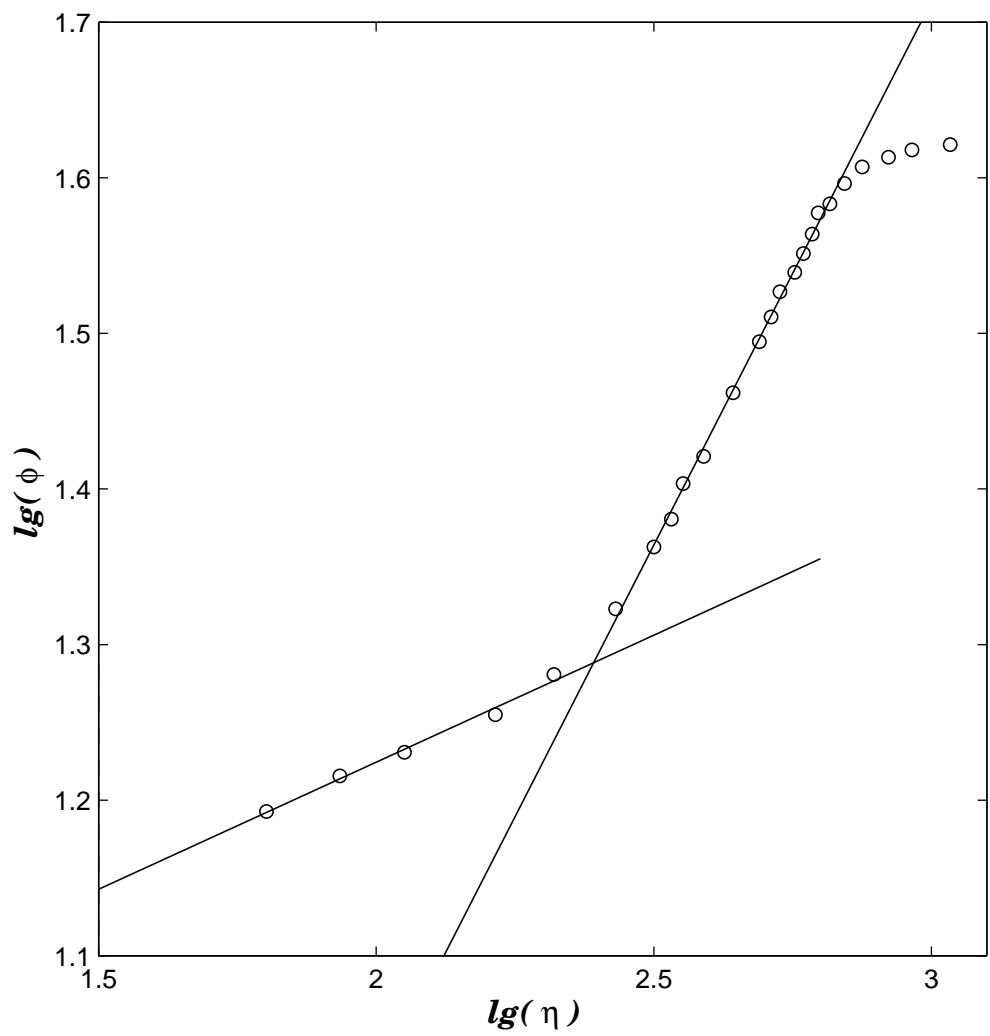
\end{document}